\documentclass[12pt]{article}
\usepackage{latexsym}
\usepackage{amsmath, amsthm}
\usepackage[all,dvips,ps]{xy}
\usepackage{amsfonts}
\usepackage{mathrsfs}
\newtheorem{prop}{Proposition}[section]
\newtheorem{lem}{Lemma}[section]
\newtheorem{thm}{Theorem}[section]
\newtheorem{cor}{Corollary}[section]

\newtheorem{exa}{Example}[section]

\newcommand{\Rs}{\mathbb{R}}

\newcommand{\N}{\mbox{null}}
\newcommand{\bz}{{\bf 0} }

\newcommand{\bpr}{{\bf Proof.} \hspace{1 em}}
\newcommand{\epr}{ \\ \hspace*{4.5in} $\Box$ }
\newcommand{\beq}{ \begin{equation} }
\newcommand{\eeq}{ \end{equation} }
\newcommand{\bt}{ \begin{tabular} }
\newcommand{\et}{ \end{tabular} }

\begin{document}

\bibliographystyle{plain}
\title{On Yielding and Jointly Yielding Entries of  Euclidean Distance Matrices }
\vspace{0.3in}
        \author{ A. Y. Alfakih
  \thanks{E-mail: alfakih@uwindsor.ca}
  \\
          Department of Mathematics and Statistics \\
          University of Windsor \\
          Windsor, Ontario N9B 3P4 \\
          Canada
}

\date{April 20, 2017. Revised \today}
\maketitle


\vspace{0.1in}

\begin{abstract}

An $n \times n$ matrix $D$ is a Euclidean distance matrix (EDM) if there exist points $p^1,\ldots,p^n$ in some Euclidean space
such that $d_{ij}=||p^i-p^j||^2$ for all $i,j=1,\ldots,n$.
Let $D$ be an EDM and let $E^{ij}$ be the $n \times n$ symmetric matrix with
1's in the $ij$th and $ji$th entries and 0's elsewhere. We say that
$[l_{ij}, u_{ij}]$ is the yielding interval of entry $d_{ij}$
if  it holds that $D+t E^{ij}$ is an EDM if and only if $l_{ij} \leq t \leq u_{ij}$.
If the yielding interval of entry $d_{ij}$ has length 0, i.e., if $l_{ij}=u_{ij}=0$, then $d_{ij}$ is said to
be unyielding. Otherwise, if $l_{ij} \neq u_{ij}$, then $d_{ij}$ is said to be yielding.
Let $d_{ij}$ and $d_{ik}$ be two unyielding entries of $D$. We say that $d_{ij}$ and $d_{ik}$ are
jointly yielding if $D+t_1 E^{ij} + t_2 E^{ik} $ is an EDM for some nonzero scalars $t_1$ and $t_2$.

In this paper, we characterize the yielding and the jointly yielding entries of an
EDM $D$ in terms of Gale transforms of  $p^1,\ldots,p^n$. Moreover,
for each yielding entry, we present explicit formulae of its yielding
interval. Finally, we specialize our results to the case where points $p^1,\ldots,p^n$ are in general position.

\end{abstract}

\section{Introduction}
An $n \times n$ matrix $D=(d_{ij})$ is said to be a Euclidean distance matrix (EDM) if there exist points $p^1,\ldots,p^n$ in some Euclidean
space such that
\[
d_{ij}= || p^i - p^j ||^2 \mbox{ for all } i,j=1,\ldots, n,
\]
where $||\;||$ denotes the Euclidean norm.
$p^1,\ldots,p^n$ are called the {\em generating points} of $D$ and the dimension of their affine span is
called the {\em embedding dimension} of $D$.
Let $D$ be a given $n \times n$ EDM and let $E^{ij}$ denote the $n \times n$
symmetric matrix with 1's in the $ij$th and $ji$th entries
and zeros elsewhere. Further, let $l_{ij} \leq 0$ and $u_{ij} \geq 0$ be the two scalars  such that
$D + t E^{ij}$ is an EDM if and only if $l_{ij} \leq t \leq u_{ij}$. That is, $D$ remains an EDM iff its $ij$th and $ji$th entries  vary
between $d_{ij}+l_{ij}$ and $d_{ij} + u_{ij}$, while keeping all other entries unchanged.
The closed interval $[l_{ij},u_{ij}]$ is called the {\em yielding interval} of entry  $d_{ij}$.
Entry $d_{ij}$ is said to be {\em unyielding} if $l_{ij}=u_{ij} = 0$;
otherwise, if $l_{ij} \neq u_{ij}$, then $d_{ij}$ is said to be {\em yielding}.
Let $d_{ij}$ and $d_{ik}$ be two unyielding entries of $D$. Then $d_{ij}$ and $d_{ik}$ are said to be
{\em jointly yielding} if $D+t_1 E^{ij}+t_2 E^{ik}$ is an EDM for some nonzero scalars $t_1$ and $t_2$.
Otherwise, they are called {\em jointly unyielding}. Note that the notion of jointly yielding (or jointly unyielding)
is defined only for two unyielding entries in the same row (column).

The theory of universal rigidity of bar frameworks provides sufficient conditions, in terms of stress matrices,
 for a given entry of $D$ to be unyielding, or for a pair of unyielding entries of $D$
 to be jointly unyielding \cite{alf16}. However, providing necessary conditions is possible only if
 further assumptions on the generating points of $D$, such as genericity, are made.
 Instead, we present in this paper a simple characterization (Theorem \ref{thmchar})
 of the unyielding entries of $D$ in terms of Gale transform of the
generating points of $D$. This characterization involves only  checking whether two vectors are parallel, without
the use of  stress matrices. Moreover, for each yielding entry of $D$, we present
explicit formulae for its yielding interval (Theorems \ref{thmir=n-1}, \ref{thmizz=0} and \ref{thmizpz}).
We also present a simple characterization (Theorem \ref{thmchar2}) of jointly unyielding entries of $D$.
This characterization involves checking whether a vector is in the linear span of two other vectors.
Finally, we specialize our results to the case where the generating points of $D$ are in general
position (Corollaries \ref{corgp1}, \ref{corgp2}, \ref{corjgp1} and \ref{corjgp2}).

\subsection{Notation}
We collect here the notation used throughout  the paper. $E^{ij}$ is the $n \times n$ matrix with 1's
 in the $ij$th and $ji$th entries and 0's elsewhere. $e$ denotes the vector of all 1's in $\Rs^n$ and $I_n$ denotes
the identity matrix of order $n$. The zero matrix or zero vector of appropriate dimension is denoted  by $\bz$. For a symmetric
matrix $A$, we mean by $A \succeq \bz \; (A \succ \bz)$ that $A$ is positive semidefinite (positive definite).
``$\backslash$" denotes the set-theoretic difference.
Finally, $\N(A)$ denotes the null space of $A$.

\section{Preliminaries}

In this section we present the necessary background concerning EDMs, Gale matrices and symmetric matrices of rank 2.

\subsection{EDMs and Gale Matrices}

Let $e$ be the vector of all 1's in $\Rs^n$ and let $V$ be an $n \times (n-1)$ matrix such that
\beq \label{defV}
V^T e = \bz \mbox{ and } V^TV=I_{n-1}.
\eeq
Thus, the orthogonal projection on $e^{\perp}$, the orthogonal complement of $e$ in $\Rs^n$, is given by
\[
 J= VV^T = I_n - \frac{ee^T}{n}.
\]
For a symmetric matrix $A$, we use $A \succeq \bz$ ( $\succ \bz$) to indicate that $A$ is positive semidefinite (positive definite).
The following is a well-known characterization of EDMs \cite{sch35,yh38, gow85,cri88}.
\begin{thm}[Schoenberg 1935 \cite{sch35}, Young and Householder 1938 \cite{yh38}]
Let $D$ be an $n \times n$ symmetric matrix whose diagonal entries are all 0's. Then $D$ is an EDM
if and only if
\[
B= -\frac{1}{2} J D J \succeq \bz,
\]
in which case, the embedding dimension of $D$ is given by rank $B$.
\end{thm}

Let $D$ be an $n \times n$ EDM of embedding dimension $r$ and let $B = -JDJ/2$ be factorized as  $B =PP^T$
where $P$ is $n \times r$. Then $p^1, \ldots, p^n$, the generating points of $D$, are given by the rows of $P$.
That is,
\beq
P = \left[ \begin{array}{c} (p^1)^T \\ \vdots \\ (p^n)^T \end{array} \right].
\eeq
Thus $P$ is called a {\em configuration matrix} of $D$. Three remarks are in order here. First,
$P$ has full column rank, i.e., rank $P = r$. Second,
$B$ is the Gram matrix of the generating points of $D$ (or the Gram matrix of $D$ for short).
Third, $P^Te=\bz$ since $Be=\bz$, thus the origin coincides with the centroid of points  $p^1, \ldots, p^n$.
This fact, is crucial for the results of this paper.

Let $F=\{ B \succeq \bz: Be=\bz\}$.
It is  well known that $F$ is a face of the cone of  positive semidefinite matrices of order $n$ \cite{wsv00}; and
that $F$ is isomorphic to the cone of positive semidefinite matrices of order $n-1$. By exploiting this fact, we  will find it more
convenient to use projected Gram matrices instead of Gram matrices to represent the generating points of an EDM $D$.
To this end, the {\em projected Gram matrix} of $D$, denoted by $X$, is defined as
\[
  X = V^T B V,
\]
where $B$ is the Gram matrix of $D$ and $V$ is as defined in (\ref{defV}). Thus,
$X =V^T PP^T V$ = $-V^T D V/2$. Moreover, $B= VXV^T$ and consequently,  rank $X$ = rank~$B$ and
$X \succeq \bz$ iff $B \succeq \bz$. Furthermore, it is easy to see
that the eigenvalues of $B$ are precisely the eigenvalues of $X$ plus one additional 0 eigenvalue.

Assume that $r \leq n-2$, and let
\beq
Z = \left[ \begin{array}{c} (z^1)^T \\ \vdots \\ (z^n)^T \end{array} \right]
\eeq
 be the $n \times (n-r-1)$ matrix whose columns form a basis of
\beq
\N(\left[ \begin{array}{c} P^T \\ e^T \end{array} \right])= \N(\left[ \begin{array}{c} B \\ e^T \end{array} \right]).
\eeq
$Z$ is called a {\em Gale matrix} of $D$ and $z^1, \ldots, z^n \in \Rs^{n-r-1}$ are called
{\em Gale transforms} of $p^1, \ldots, p^n$. In fact, the columns of $Z$ express the affine dependency of the
generating points of $D$.
Gale transform \cite{gal56, gru67} is well known and widely used in the theory of polytopes.
\begin{exa} \label{ex1}
Consider the EDM  $D = \left[ \begin{array}{ccccc} 0 & 1 & 4 & 2 & 2 \\ 1 & 0 & 1 & 1 & 1 \\
                                                   4 & 1 & 0 & 2 & 2 \\ 2 & 1 & 2 & 0 & 4 \\
                                                   2 & 1 & 2 & 4 & 0 \end{array} \right] $ of embedding dimension~$2$.
A configuration matrix of $D$ is $P=\left[ \begin{array}{rr} -1 & 0 \\ 0 & 0 \\ 1 & 0 \\ 0 & 1 \\ 0 & -1 \end{array} \right]$.
Thus a Gale matrix of $D$ is $Z=\left[ \begin{array}{rr} 1 & 0 \\ -2 & -2 \\ 1 & 0 \\ 0 & 1 \\ 0 & 1 \end{array} \right]$.
Therefore, the generating points of $D$ are
\[ p^1= \left[ \begin{array}{r} -1 \\ 0 \end{array} \right],
p^2= \left[ \begin{array}{r} 0 \\ 0 \end{array} \right],
p^3 = \left[ \begin{array}{r} 1 \\ 0 \end{array} \right],
p^4= \left[ \begin{array}{r} 0 \\ 1 \end{array} \right] \mbox{ and }
p^5= \left[ \begin{array}{r} 0 \\ -1 \end{array} \right],
\]
and their Gale transforms are
\[ z^1= \left[ \begin{array}{r} 1 \\ 0 \end{array} \right],
z^2= \left[ \begin{array}{r} -2 \\ -2 \end{array} \right],
z^3 = \left[ \begin{array}{r} 1 \\ 0 \end{array} \right],
z^4= \left[ \begin{array}{r} 0 \\ 1 \end{array} \right] \mbox{ and }
z^5= \left[ \begin{array}{r} 0 \\ 1 \end{array} \right].
\]
\end{exa}

The following lemma, relating Gale matrices and the projected Gram matrix, is crucial for the proofs of our results.
\begin{lem}[\cite{alf07}] \label{lemVUZ}
Let $X$ be the projected Gram matrix of an $n \times n$ EDM  $D$ of embedding dimension $r$. Let $Z$ and $P$ be, respectively, a Gale matrix and
a configuration matrix of $D$, where $P^Te=\bz$. Further,
let $U$ and $W$ be the matrices whose columns form orthonormal bases of the null space and the column space of $X$, respectively.
Then
\begin{enumerate}
\item $ VU = ZA$, where $A$ is an $(n-r-1) \times (n-r-1)$ nonsingular matrix,
\item $ VW = PA'$, where $A'$ is an $r \times r$ nonsingular matrix,
\end{enumerate}
where $V$ is as defined in (\ref{defV}).
\end{lem}

To illustrate part 2 of Lemma \ref{lemVUZ}, let
$X = W \Lambda W^T$ be the spectral decomposition of $X$, where $\Lambda$ is the $r \times r$ diagonal matrix
consisting of the positive eigenvalues of $X$. Thus $P = V W \Lambda^{1/2}$ is a configuration matrix of $D$.
Hence, $P^TP= \Lambda$ and $A'=\Lambda^{-1/2}$ in this case.

The following proposition summarizes few useful properties of Gale transform. It plays a crucial role in establishing
the yielding intervals of the entries of an EDM $D$.

\begin{prop} \label{propz=0}
Let $D$ be an $n \times n$, $n \geq 3$,  EDM of embedding dimension $r \leq n-2$.
Let $Z$ and $P$ be, respectively, a Gale matrix and a configuration matrix of $D$,
where $P^Te=\bz$. Let $i$, $j$ and $k$ be three distinct indices in $\{1,\ldots,n\}$.
\begin{enumerate}
\item If $z^i=\bz$, then $p^i \neq \bz$ and $p^i$ is not in the affine hull of $\{p^1, \ldots, p^n \} \backslash \{p^i\}$.
\item If $n \geq 4$ and if $z^i= z^j=\bz$, then $p^i \neq \bz $, $p^j \neq \bz$ and $p^i \neq c' p^j$ for any  scalar $c'$.
\item If $z^i \neq \bz$, $z^j \neq \bz$ and $z^i= c z^j$ for some nonzero scalar $c$, then  $p^i - c p^j \neq \bz$.
\item If $n \geq 4$ and if $z^i \neq \bz$, $z^j \neq \bz$, $z^k \neq \bz$ and $z^i= c_1 z^j + c_2 z^k$ for some
nonzero scalars $c_1$ and $c_2$, then  $p^i - c_1 p^j - c_2 p^k \neq \bz$.
\end{enumerate}
\end{prop}

\bpr

Wlog assume that $z^1=0$ and assume, to the contrary, that $p^1$ is in the affine hull of $\{p^2,\ldots,p^n\}$. Then
there exist scalars $\lambda_2, \ldots,\lambda_n$ such that
\[
\left[ \begin{array}{c} p^1 \\ 1 \end{array} \right] = \sum_{i=2}^n
\lambda_i \left[ \begin{array}{c} p^i \\ 1 \end{array} \right].
\]
Let $x=[-1  \;\; \lambda_2 \;\; \cdots \;\; \lambda_n]^T$ in $\Rs^n$.
Thus, there exists $\xi \neq \bz$ in $\Rs^{n-r-1}$ such that $Z\xi = x$. Hence, $(z^1)^T \xi = -1$, a contradiction.
Therefore, $p^1$ is not in the affine hull of $\{p^2,\ldots,p^n\}$.
Now assume, to the contrary, that $p^1=\bz$. Thus  $p^1 = \bz = p^2+ \cdots +p^n$
since $P^Te=\bz$. Hence,  $p^1$ is in the affine hull of $\{p^2,\ldots,p^n\}$, a contradiction.
Therefore, $p^1 \neq \bz$, and  part 1 is proven.

To prove part 2, wlog assume that $z^1=z^2=\bz$ and assume, to the contrary, that $p^1 = c'p^2$ for some scalar $c'$.
Then, it follows from part 1 that $p^1 \neq \bz$ and $p^2 \neq \bz$ and hence, $c' \neq 0$.
Since $P^Te=\bz$, it follows that
$p^1 = c' p^2 + \beta \sum_{i=1}^n p^i$, where $\beta$ is an arbitrary  scalar.
Hence, $(1-\beta)p^1 = (c' + \beta) p^2 + \beta \sum_{i=3}^n p^i$.
Let $x=[-1+\beta  \;\; c'+\beta \;\; \beta \;\; \cdots \;\; \beta]^T$ in $\Rs^n$.
If we set  $\beta = (1-c')/n$, then $x \neq \bz$ and
$e^T x = 0$. Hence, there
exists $\xi \neq \bz$  such that $Z\xi = x$, and in particular, $(z^1)^T \xi = 0 = -1+\beta$ and $(z^2)^T \xi = 0 = c'+\beta$.
Thus $-c'=\beta=1$ and hence, $n=2$, a contradiction. Therefore,
there does not exist a scalar $c'$ such that $p^1 = c' p^2$.

To prove part 3, wlog assume that $z^1 \neq \bz$, $z^2 \neq \bz$,  $z^1=cz^2$, $c \neq 0$
and assume, to the contrary, that $p^1=cp^2$.
Let $x=[-1 \;\; c \;\; 0 \;\; \cdots \;\; 0]^T$ in $\Rs^n$. Then $Z^T x = \bz$ and $P^T x = \bz$.
But, $\N( \left[ \begin{array}{c} Z^T \\ P^T \end{array} \right] )$ = span of $e$ since $P^Te=\bz$.
Hence,
$x= \alpha e$ for some scalar $\alpha$, a contradiction since $x$ has at least one zero entry. Therefore, $p^1 - c p^2 \neq \bz$.

The proof of part 4 is similar to that of part 3 where in this case
$x=[-1 \;\; c_1 \;\; c_2 \;\; 0 \;\; \cdots \;\; 0]^T$
 \epr

We remark here that in part 3 of Proposition \ref{propz=0},
$p^i$ may be parallel to $p^j$, say $p^i=c'p^j$, but
$c'$ cannot be equal to $c$ as illustrated by the following example.

\begin{exa} \label{ex2}
Consider the EDM  $D = \left[ \begin{array}{cccc} 0 & 0 & 1 & 1 \\ 0 & 0 & 1 & 1 \\
                                                   1 & 1 & 0 & 4  \\ 1 & 1 & 4 & 0  \end{array} \right] $ of embedding dimension~$1$.
A configuration matrix of $D$ is $P=\left[ \begin{array}{r}  0 \\ 0  \\ - 1 \\  1  \end{array} \right]$.
Thus a Gale matrix of $D$ is $Z=\left[ \begin{array}{rr} -2 & 0 \\ 0 & -2 \\ 1 & 1 \\ 1 & 1 \end{array} \right]$.
Note that $z^3=z^4$ and $p^3 = - p^4$; i.e, $p^3$ is parallel to $p^4$ but $c' \neq c$.
\end{exa}

\subsection{A Property of Symmetric Rank-Two Matrices}

Vectors $u$ and $v$ in $\Rs^n$ are {\em parallel} if $u = c v$ for some nonzero scalar $c$. Thus,
if $u=v=\bz$, then $u$ and $v$ are parallel. The following proposition will be quite useful in the sequel.

\begin{prop} \label{prop1}
Let $a$ and $b$ be two nonzero, nonparallel vectors in $\Rs^r$, $r \geq 2$, and let $\Psi=ab^T+ba^T$. Then $\Psi$ has
exactly one positive eigenvalue $\lambda_1$ and one negative eigenvalue $\lambda_r$, where
\[
\lambda_1 = a^Tb + ||a|| \; ||b||   \mbox{ and } \lambda_r = a^Tb - ||a|| \; ||b||.
\]
\end{prop}

\bpr
Assume that $r=2$ and let the eigenvectors of $\Psi$ be of the form  $xa+yb$, where $x$ and $y$ are scalars.
Then $\Psi (xa+yb) = \lambda (xa+yb)$ leads to the following system of equations
\beq \label{eqrank2}
\left[ \begin{array}{cc} a^Tb & ||b||^2  \\ ||a||^2 & a^Tb \end{array} \right]
\left[ \begin{array}{c} x \\ y \end{array} \right] = \lambda
\left[ \begin{array}{c} x \\ y \end{array} \right].
\eeq
Hence, the eigenvalues of $\Psi$ are precisely the eigenvalues of
$\left[ \begin{array}{cc} a^Tb & ||b||^2  \\ ||a||^2 & a^Tb \end{array} \right]$, which are
$\lambda_1 = a^Tb + ||a|| \; ||b||$   and $\lambda_r = a^Tb - ||a|| \; ||b||$.

Now assume that $r \geq 3$ and let $u^1, \ldots, u^{r-2}$ be an orthonormal basis of the null space of
$\left[ \begin{array}{c} a^T \\ b^T \end{array} \right]$. Then obviously, $u^1, \ldots, u^{r-2}$ are orthonormal eigenvectors
of $\Psi$ corresponding to eigenvalue 0. Thus we have 2 remaining eigenvectors of $\Psi$ of the form $xa+yb$, where
$x$ and $y$ satisfy Equation (\ref{eqrank2}). Therefore, the remaining two eiegvalues of $\Psi$ are $\lambda_1$ and $\lambda_r$
as given above.
The fact that $\lambda_1 > 0 $ and $\lambda_r < 0$
follows from the Cauchy-Schwarz inequality since $a$ and $b$ are nonzero and nonparallel.
\epr

\section{Characterizing the Unyielding Entries}

We consider, first, the case where the generating points of $D$ are affinely independent.

\begin{thm}\label{thmr=n-1}
Let $D$ be an $n \times n$ EDM of embedding dimension $r = n-1$. Then every entry of $D$ is yielding.
\end{thm}

\bpr
Let $1 \leq k< l \leq n$. Then
$D + t E^{kl}$ is EDM iff $2X - t V^T E^{kl} V \succeq \bz$, where $X$ is the projected Gram matrix of $D$
and $V$ is as defined in (\ref{defV}).
But, $X \succ \bz$  since $X$ is of order $n-1$ and rank $X$ = $r=n-1$. Thus, obviously,
there exists $t \neq 0$ such that $2X - t V^T E^{kl} V \succeq \bz$.
Consequently, $d_{kl}$ is yielding and the result follows.
\epr

The following lemma is needed for the case where the embedding dimension of $D$ is $r \leq n-2$.

\begin{lem} \label{lembasic}
Let $D$ be an $n \times n$ nonzero EDM of embedding dimension $r \leq n-2$, and let $Z$ and $P$ be a Gale matrix and
a configuration matrix of $D$, respectively, where $P^Te=\bz$. Further, let $X$ be the projected Gram matrix $D$. Then
$2X - t V^T E^{kl} V \succeq \bz$ iff
\[
\left[ \begin{array}{rr} 2(P^TP)^2 - t \;( p^k (p^l)^T + p^l (p^k)^T) &
                         -t \; (p^k (z^l)^T + p^l (z^k)^T) \\ -t \; (z^k (p^l)^T + z^l (p^k)^T) &
                         - t \; (z^k (z^l)^T + z^l (z^k)^T) \end{array} \right] \succeq \bz.
\]

\end{lem}

\bpr
Let $W$ and $U$ be the two matrices whose columns form orthonormal bases of the column space and the null space of $X$,
respectively, and thus $Q=[W \;\; U]$ is orthogonal. Hence,
$2X - t V^T E^{kl} V \succeq \bz$ iff
\[
Q^T (2 X - t V^T E^{kl}V ) Q = \left[ \begin{array}{rr} 2 W^T X W  - t W^T V^T E^{kl} VW &
                   -t W^T V^T E^{kl} V U \\ -t U^T V^T E^{kl} V W & -t U^T V^T E^{kl}V U \end{array} \right] \succeq \bz.
\]

But it follows from Lemma \ref{lemVUZ} that $VU = Z A $ and $VW=P A'$, where $A$ and $A'$ are nonsingular. Hence,
$2X - t V^T E^{kl} V \succeq \bz$ iff
\[
\left[ \begin{array}{rr} 2(P^TP)^2 - t \; (p^k (p^l)^T + p^l (p^k)^T) &
                   -t \;( p^k (z^l)^T + p^l (z^k)^T) \\ -t \; (z^k (p^l)^T + z^l (p^k)^T) &
                   - t \; (z^k (z^l)^T + z^l (z^k)^T) \end{array} \right] \succeq \bz.
\]
\epr

It is worth pointing out, here, that $(P^TP)^2 \succ \bz$.
Now we are ready to present and prove the theorem characterizing the unyielding entries of $D$.

\begin{thm} \label{thmchar}
Let $D$ be an $n \times n$ EDM of embedding dimension $r \leq n-2$, and let $z^1,\ldots, z^n$ be Gale transforms of the generating points of $D$.
Then entry $d_{kl}$ is unyielding if and only if $z^k$ is not parallel to $z^l$; i.e., iff
there does not exist a nonzero scalar $c$ such that $z^k = c z^l$.
\end{thm}

\bpr
Let $1 \leq k<l \leq n$. Entry
$d_{kl}$ is yielding iff there exists $t\neq 0$ such that $D+t E^{kl}$ is an EDM or equivalently, iff $2X - t V^T E^{kl}V \succeq \bz$,
where $X$ is the projected Gram matrix of $D$, and $V$ is as defined in (\ref{defV}).

Assume that $z^k$ is parallel to $z^l$, i.e.,  $z^k = c z^l$ for some nonzero scalar $c$.
Then $z^k (z^l)^T + z^l (z^k)^T = 2 c z^l (z^l)^T \succeq 0 $ and $p^k (z^l)^T + p^l (z^k)^T = (p^k + c p^l) (z^l)^T$. Hence,
$\N( z^l (z^l)^T)$ = $\N((z^l)^T) \subseteq $ $\N( (p^k+c p^l) (z^l)^T)$.
Therefore, it follows from Lemma~\ref{lembasic} that there exists $t \neq  0$ such that
$2X - t V^T E^{kl} V \succeq \bz$ and thus $d_{kl}$ is yielding.

On the other hand, assume that $z^k$ and $z^l$ are not parallel and assume, to the contrary, that entry $d_{kl}$ is yielding.
Then there exists $t \neq 0$ such that $2X - t V^T E^{kl} V \succeq \bz$. Thus it follows from Lemma \ref{lembasic} that
there exists $t \neq 0$ such that $-t (z^k (z^l)^T + z^l (z^k)^T) \succeq \bz$ and
$\N( z^k (z^l)^T+z^l (z^k)^T) \subseteq $ $\N( p^k (z^l)^T + p^l (z^k)^T)$.
Next, we consider two cases:

(i)  $z^k=\bz$ and $z^l \neq \bz$. In this case,
$z^k (z^l)^T + z^l (z^k)^T = \bz$ and  $p^k (z^l)^T + p^l (z^k)^T) = p^k (z^l)^T \neq \bz$ since,
by Proposition \ref{propz=0} (part 1), $p^k \neq \bz$. Hence, we have a contradiction since $\N(\bz) \not \subseteq \N(p^k (z^l)^T)$.

(ii) both $z^k$ and $z^l$ are nonzero. Also, in this case we have a contradiction since Proposition \ref{prop1} implies
that $z^k (z^l)^T+z^l (z^k)^T$ is indefinite.
Therefore, $d_{kl}$ is unyielding.
\epr

\begin{exa} \label{ex11}
Let $D$ be the EDM considered in Example \ref{ex1}. Then
$d_{13}$ and $d_{45}$ are yielding, while all other entries are unyielding.
\end{exa}

Points $p^1,\ldots, p^n$ in $\Rs^r$ are said to be in {\em general position} if every $r+1$ of them are affinely independent.
For instance, points in the plane are in general position if no three of them are collinear. An EDM $D$ of embedding dimension $r$
is said to be in general position if its generating points are in general position in $\Rs^r$.
Points in general position have a nice characterization in terms of Gale matrices.

\begin{lem}[\cite{alf07}] \label{lemgp}
Let $D$ be an $n \times n$ EDM of embedding dimension $r$, $r \leq n-2$. Let $Z$ be a Gale matrix of $D$.
Then $D$ is in general position if and only if every submatrix of $Z$ of order $(n-r-1)$ is nonsingular.
\end{lem}

\begin{cor} \label{corgp1}
Let $D$ be an $n \times n$ EDM of embedding dimension $r=n-2$. Then $D$ is in general position in $\Rs^r$ if
and only if every entry of $D$ is yielding.
\end{cor}

\bpr
In this case, $z^1, \ldots, z^n$ are scalars since $(n-r-1)=1$. Assume that $D$ is in general position.
Then it follows from Lemma \ref{lemgp} that
$z^1, \ldots, z^n$ are nonzero. Thus $z^k$ is parallel to $z^l$
for all $1 \leq k < l \leq n$, and hence every entry of $D$ is yielding.

Now assume that one entry of $D$ say, $d_{kl}$, is unyielding. Then $z^k$ is not parallel to $z^l$. Thus
either $z^k=0$ or $z^l=0$ but not both. Therefore,  it follows from Lemma \ref{lemgp} that $D$ is not in general position.
\epr

\begin{cor} \label{corgp2}
Let $D$ be an $n \times n$ EDM of embedding dimension $r\leq n-3$. If $D$ is in general position in $\Rs^r$,
then every  entry of $D$ is unyielding.
\end{cor}

\bpr
Assume, to the contrary, that one entry of $D$, say $d_{kl}$, is yielding.
Thus, it follows from Theorem \ref{thmchar} that $z^k=cz^l$ for some nonzero scalar $c$. Note that in this
case, $(n-r-1) \geq 2$.
Hence, any $(n-r-1) \times (n-r-1)$ submatrix of $Z$ containing $(z^k)^T$ and $(z^l)^T$ is singular. This
contradicts  Lemma \ref{lemgp} and the proof is complete.
\epr

The following example shows that the converse of Corollary \ref{corgp2} is not true.

\begin{exa}\label{exalast}
Consider the EDM  $D = \left[ \begin{array}{ccccc} 0 & 1 & 4 & 9 & 1 \\ 1 & 0 & 1 & 4 & 0 \\
                                                   4 & 1 & 0 & 1 & 1  \\ 9 & 4 & 1 & 0 & 4 \\
                                                   1 & 0 & 1 & 4 & 0 \end{array} \right] $ of embedding dimension~$1$.
A configuration matrix and a Gale matrix  of $D$ are
\[ P=\frac{1}{5} \left[ \begin{array}{r} -7 \\ -2 \\ 3 \\ 8 \\ -2 \end{array} \right] \mbox{ and }
 Z= \left[ \begin{array}{rrr} 1 & 0 & 0 \\ -2 & 1 & -1 \\ 1 & -2 & 0 \\  0 & 1 & 0 \\ 0 & 0 & 1 \end{array} \right].
\]
Obviously, $D$ is not in general position in $\Rs^1$ since $p^2=p^5$. However, every entry of $D$ is unyielding.
\end{exa}

Finally, an important consequence of Theorem \ref{thmr=n-1} and Corollaries \ref{corgp1} and \ref{corgp2}
is worth pointing out. If an $n \times n$ EDM $D$ of embedding dimension $r$ is in general position,
 then the entries of $D$ are either all
yielding (if  $n=r+1$ or $n=r+2$) or all unyielding (if $n \geq r+3$).

\section{Determining the Yielding Intervals}

Let $d_{kl}$ be a yielding entry of $D$. To determine the yielding interval of $d_{kl}$, we treat
each of the following three cases separately:
First, when $r=n-1$. Second, when $r \leq n-2$ and $z^k=z^l=\bz$.
Third, when $r \leq n-2$, $z^k \neq \bz$, $z^l \neq \bz$ and $z^k=cz^l$ for some nonzero scalar $c$. As will be seen below, the first two cases
are similar. More specifically, in the first two cases, $0$ is in the interior of the yielding interval, while in the third case,
$0$ is an endpoint.

Let $D$ be an EDM of order 2. Then the off-diagonal entry of $D$ can assume any nonnegative value.
As a result, the case of EDMs of order 2 is trivial and uninteresting. Consequently,
we focus in what follows on EDMs of order $\geq 3$.

\begin{prop}\label{prop3}
Let $D$ be  an $n \times n $, $n \geq 3$, EDM matrix of embedding dimension $r=n-1$ and let
$P$ be a configuration matrix of $D$, where $P^Te=\bz$. Let $S=P (P^TP)^{-1}$ and let $(s^i)^T$ be the $i$th row of $S$; i.e.,
$s^i = (P^TP)^{-1} p^i$.
Then $s^k$ and $s^l$ are nonzero and nonparallel for all $k \neq l$.
\end{prop}

\bpr
Assume, to the contrary, that $s^k=\bz$. Then $p^k=\bz$ and thus
$P^T e^k = \bz$, where $e^k$ is the $k$th standard unit vector in $\Rs^n$.
Since $P^Te=\bz$, this
implies that rank ($PP^T$) = $r \leq n-2$, a contradiction. To complete the proof, assume, to the contrary,
that $s^k = c s^l$ for some nonzero scalar $c$, where $k \neq l$. Then $p^k = c p^l$ and thus
$P^T (e^k-c e^l)=\bz$ and again we have a contradiction.
 \epr

\begin{thm} \label{thmir=n-1}
Let $D$ be an $n \times n$, $n \geq 3$, EDM of embedding dimension $r=n-1$ and let $P$  be a configuration matrix a
of $D$, where $P^Te=\bz$. Further, Let $S= P (P^TP)^{-1}$ and let  $(s^i)^T$ be the $i$th row of $S$.
Then the yielding interval of entry $d_{kl}$ is given by
\[
\left[ \frac{2}{\lambda_{r}} \; , \; \frac{2}{\lambda_1} \right],
\]
where $\lambda_1 = (s^k)^T s^l + ||s^k|| \; ||s^l||$ and
$\lambda_r = (s^k)^T s^l - ||s^k|| \; ||s^l||$.
\end{thm}

\bpr Let $1 \leq k < l \leq n$ and let $X$ be the projected Gram matrix of $D$.
Let $X = W \Lambda W^T$ be the spectral decomposition of $X$. Thus $\Lambda \succ \bz$ and $W$ is orthogonal since $r=n-1$.
Thus $D + t E^{kl}$ is EDM iff $2X - t V^T E^{kl} V \succeq \bz$ iff $2 W^T X W - t W^T V^T E^{kl} V W \succeq \bz$.
But $W^TXW = W^T V^T PP^T V W$. Thus, it follows from Lemma \ref{lemVUZ} that
$2 W^T X W - t W^T V^T E^{kl} V W \succeq \bz$ iff
\beq \label{eq1}
2 (P^T P )^2  - t P^T E^{kl} P \succeq \bz.
\eeq
But Equation (\ref{eq1}) holds iff
\[
2I_{n-1}  - t (P^TP)^{-1} P^T E^{kl} P  (P^TP)^{-1} = 2 I_{n-1} - t S^T E^{kl} S \succeq \bz.
\]
In light of Propositions \ref{prop1} and \ref{prop3},
let $\lambda_1 > 0 > \lambda_r$ be the nonzero eigenvalues of $S^T E^{kl} S$ = $s^k (s^l)^T + s^l (s^k)^T$. Thus,
$2 I_{n-1} - t S^T E^{kl} S \succeq \bz$ iff $2-t\lambda_1 \geq 0$ and $2-t\lambda_{r} \geq 0$.
Therefore, $D + t E^{kl}$ is EDM iff $2/\lambda_{r} \leq t \leq 2/\lambda_1$.
\epr

Note that if $P$ and $P'$ are two configuration matrices of $D$ such that $P^Te={P'}^Te=\bz$. Then $P'=PQ$ for some
orthogonal matrix $Q$. Thus $S'= P' ({P'}^T P')^{-1}$ = $P ({P}^T P)^{-1}Q$ = $S Q$. Thus
${S'}^T E^{kl} S' = Q^T{S}^T E^{kl} S Q $ and hence, the matrices ${S'}^T E^{kl} S'$ and ${S}^T E^{kl} S$ are similar.

\begin{exa}
Consider the EDM  $D = \left[ \begin{array}{ccc} 0 & 4 & 10 \\ 4 & 0 & 10 \\
                                                   10 & 10 & 0  \end{array} \right] $ of embedding dimension~$2$.
A configuration matrix of $D$ is $P=\left[ \begin{array}{rr}  -1 & -1 \\ 1 & - 1 \\  0 & 2  \end{array} \right]$.
Thus $S=P (P^TP)^{-1} = \left[ \begin{array}{rr}  -1/2 & -1/6 \\ 1/2 & - 1/6 \\  0 & 1/3  \end{array} \right]$.
Hence, $S^T E^{12} S = \left[ \begin{array}{cc}  -1/2 & 0 \\ 0 & 1/18  \end{array} \right]$ with nonzero eigenvalues
$ (s^1)^T s^2 - ||s^1|| \; ||s^2|| = -1/2$ and $(s^1)^T s^2 + ||s^1|| \; ||s^2|| = 1/18$,
$S^T E^{13} S = \left[ \begin{array}{cc}  0 & -1/6 \\ -1/6 & - 1/9  \end{array} \right]$
with nonzero eigenvalues $(s^1)^T s^3 - ||s^1|| \; ||s^3||=(-1-\sqrt{10})/18$ and
$ (s^1)^T s^3 + ||s^1|| \; ||s^3||= (\sqrt{10}-1)/18$; and
$S^T E^{23} S = \left[ \begin{array}{cc}  0  & 1/6 \\ 1/6 & -1/9 \end{array} \right]$ with nonzero eigenvalues
$(s^3)^T s^2 - ||s^3|| \; ||s^2||=(-1-\sqrt{10})/18$ and
$ (s^3)^T s^2 + ||s^3|| \; ||s^2||=(\sqrt{10}-1)/18$. Therefore, the yielding interval of $d_{12}$
is $[-4 \;, \; 36]$, while the yielding interval of $d_{13}$ and $d_{23}$ is
\[
\left[ \frac{-36}{\sqrt{10}+1} \; , \; \frac{36}{\sqrt{10}-1} \right] = \left[ -4 \sqrt{10} + 4  \; , \; 4\sqrt{10} + 4 \right].
\]
Note that, in this case, these intervals could have been calculated using triangular inequalities.
\end{exa}

\begin{thm}\label{thmizz=0}
Let $D$ be an $n \times n$, $n \geq 4$, EDM of embedding dimension $r \leq n-2$ and let $Z$ and $P$  be a
Gale matrix and a configuration matrix of $D$ respectively, where $P^Te=\bz$.
Further, Let $S= P (P^TP)^{-1}$ and let  $(s^i)^T$ be the $i$th row of $S$.
If $z^k=z^l=\bz$, then
the yielding interval of entry $d_{kl}$ is given by
\[
\left[ \frac{2}{\lambda_{r}} \; , \; \frac{2}{\lambda_1} \right],
\]
where $\lambda_1 = (s^k)^T s^l + ||s^k|| \; ||s^l||$ and
$\lambda_r = (s^k)^T s^l - ||s^k|| \; ||s^l||$.
\end{thm}

\bpr It follows from Proposition \ref{propz=0} (part 2) that
$s^k$ and $s^l$ are nonzero and nonparallel. Moreover,
in this case
\[
\left[ \begin{array}{rr} 2(P^TP)^2 - t \; (p^k (p^l)^T + p^l (p^k)^T) &
                   -t \;( p^k (z^l)^T + p^l (z^k)^T) \\ -t \; (z^k (p^l)^T + z^l (p^k)^T) &
                   - t \; (z^k (z^l)^T + z^l (z^k)^T) \end{array} \right]
\]
reduces to
\[
\left[ \begin{array}{cr} 2(P^TP)^2 - t \; P^T E^{kl} P & \bz \\ \bz & \bz \end{array} \right].
\]
Using Proposition \ref{prop1} and Lemma \ref{lembasic}, the proof proceeds along the same lines of the proof of Theorem~\ref{thmir=n-1}.
\epr

\begin{exa}
Consider the EDM  $D = \left[ \begin{array}{cccc} 0 & 0 & 1 & 1 \\ 0 & 0 & 1 & 1 \\
                                                   1 & 1 & 0  & 2 \\ 1 & 1 & 2 & 0 \end{array} \right] $ of embedding dimension~$2$.
A configuration matrix of $D$ is $P=\frac{1}{4}\left[ \begin{array}{rr}  -1 & -1 \\ -1 & -1  \\ 3 & -1 \\ -1 & 3  \end{array} \right]$
and a Gale matrix of $D$ is $Z= \left[ \begin{array}{r}  -1  \\ 1  \\ 0 \\ 0  \end{array} \right]$.
Thus, entries $d_{12}$ and $d_{34}$ are yielding, while all other entries are unyielding. Notice that $d_{12}$ does not
fall into this case since $z^1$ and $z^2$ are nonzero. Thus we consider, next, $d_{34}$ only.

 $S=P (P^TP)^{-1} = \left[ \begin{array}{rr}  -1/2 & -1/2 \\ -1/2 & - 1/2 \\  1 & 0 \\ 0 & 1  \end{array} \right]$.
Hence, $S^T E^{34} S = \left[ \begin{array}{cc}  0 & 1 \\ 1 & 0  \end{array} \right]$ with nonzero eigenvalues
$(s^3)^T s^4 - ||s^3|| \; ||s^4||=-1$ and $(s^3)^T s^4 + ||s^3|| \; ||s^4||=1$.
Therefore, the yielding interval of $d_{34}$ is $[-2 \;, \; 2]$.
\end{exa}

\begin{thm}\label{thmizpz}
Let $D$ be an $n \times n$, $n \geq 3$, EDM of embedding dimension $r \leq n-2$ and let $Z$ and $P$ be a
Gale matrix and a configuration matrix of $D$ respectively, where $P^Te=\bz$.
Further, Let $S= P (P^TP)^{-1}$ and let $s^i$ be the $i$th row of $S$, i.e., $s^i = (P^TP)^{-1} p^i$.
If both $z^k$ and $z^l$ are nonzero and $z^k=cz^l$ for some nonzero scalar $c$ , then
the yielding interval of entry $d_{kl}$ is given by
\[
\left[\frac{- 4 c}{|| s^k - c s^l||^2} \; , \; 0 \right]  \mbox{ if } c > 0,
\]
and
\[
\left[0 \; , \; \frac{4 \; |c|}{|| s^k - c s^l||^2} \right]  \mbox{ if } c < 0.
\]
\end{thm}

\bpr
Assume that $z^k$ and $z^l$ are nonzero and $z^k=cz^l$, where $c \neq 0$.
 Let $1 \leq k < l \leq n$ and let $X$ be the projected Gram matrix of $D$.
 The $D + t E^{kl}$ is an EDM iff $2X - t V^T E^{kl} V \succeq \bz$. In light of Lemma \ref{lembasic},
 $D + t E^{kl}$ is an EDM iff
\beq \label{eq2}
 \left[ \begin{array}{cc} 2(P^TP)^2 - t \; (p^k (p^l)^T + p^l (p^k)^T) &
                   -t \;( p^k + c p^l) (z^l)^T) \\ -t \; (z^l (p^k+c p^l)^T) &
                   - t \; 2 c z^l (z^l)^T \end{array} \right] \succeq \bz.
\eeq
Assume that $r = n-2$, i.e., $(n-r-1)=1$. Then $z^l$ is a nonzero scalar. Using
Schur Complement, we have that Equation (\ref{eq2}) holds iff
\beq \label{eqnadd}
tc \leq 0 \mbox{ and } 2(P^TP)^2 + \frac{t}{2c} (p^k - c p^l)(p^k-cp^l)^T \succeq \bz,
\eeq
which is equivalent to
\[
tc \leq 0 \mbox{ and } 2I_r  + \frac{t}{2c} (s^k - c s^l)(s^k-cs^l)^T \succeq \bz,
\]
which in turn is equivalent to
\[
tc \leq 0 \mbox{ and } 2 + \frac{t}{2c} || s^k-cs^l||^2 \geq 0.
\]
The result follows from Proposition \ref{propz=0} (part 3) since $s^k - c s^l \neq \bz$.

Now assume that $r \leq n-3$, i.e., $(n-r-1) \geq 2$.
Let $Q'=[ \frac{z^l}{||z^l||} \; \;  M]$ be an $ (n-r-1) \times (n-r-1)$
orthogonal matrix and define the $ (n-1) \times (n-1)$ matrix
$Q= \left[ \begin{array}{cc} I_r & \bz \\ \bz & Q' \end{array} \right]$. Thus obviously $Q'$ is orthogonal.
By multiplying the LHS of Equation (\ref{eq2}) from the left with $Q^T$ and from the right with $Q$ we get that
 $D + t E^{kl}$ is an EDM iff
\beq \label{eq3}
 \left[ \begin{array}{cc} 2(P^TP)^2 - t \; (p^k (p^l)^T + p^l (p^k)^T) &
                   -t \;( p^k + c p^l) \; ||z^l||) \\ -t \; ||z^l|| \; (p^k+c p^l)^T) &
                   - t \; 2 c \; ||z^l||^2 \end{array} \right] \succeq \bz.
\eeq
Again using Schur complement we arrive at Equation (\ref{eqnadd}) and thus the proof is complete.
\epr

\begin{exa}
Let $D$ be the EDM considered in Example \ref{ex1}. Then
$S = P (P^TP)^{-1} = \frac{1}{2} P$.
For yielding entry $d_{13}$ we have $z^1 =  z^3$, thus $c=1$ and $s^1-cs^3$ =
$\left[ \begin{array}{r}-1  \\ 0  \end{array} \right]$. Thus the yielding interval for $d_{13}$ is $[-4 \; , \; 0]$.
On the other hand, for yielding entry $d_{45}$ we have $z^4 =  z^5$, thus $c=1$ and $s^4-cs^5$ =
$\left[ \begin{array}{r}0  \\ 1 \end{array} \right]$. Thus, the yielding interval for $d_{45}$ is also  $[-4 \; , \; 0]$.
\end{exa}

\section{Jointly Yielding Entries}

The following theorem characterizes the jointly yielding entries of an EDM $D$.

\begin{thm} \label{thmchar2}
Let $D$ be an $n \times n$, $n \geq 4$, EDM of embedding dimension $r \leq n-2$, and
let $z^1,\ldots, z^n$ be Gale transforms of the generating points of $D$.
Further, let $d_{ij}$ and $d_{ik}$ be two unyielding entries of $D$. Then
$d_{ij}$ and $d_{ik}$ are jointly yielding  if and only if there exist nonzero scalars $c_1$ and $c_2$ such that
$z^i= c_1z^j + c_2 z^k$.
\end{thm}

\bpr
The proof is analogous to that of Theorem \ref{thmchar}.
Let $d_{ij}$ and $d_{ik}$ be unyielding entries of $D$. Then $d_{ij}$ and $d_{ik}$ are jointly yielding
iff there exist $t_1\neq 0$ and $t_2 \neq 0$ such that
$D+t_1 E^{ij}+ t_2 E^{ik}$ is an EDM or equivalently, iff
\beq \label{eqj1}
\left[ \begin{array}{rr} 2(P^TP)^2 - P^T(t_1 E^{ij}+ t_2 E^{ik})P & - P^T(t_1 E^{ij}+ t_2 E^{ik})Z \\
                                      - Z^T(t_1 E^{ij}+ t_2 E^{ik})P   & - Z^T(t_1 E^{ij}+ t_2 E^{ik})Z
                                        \end{array} \right] \succeq \bz.
\eeq
Assume that $z^i = c_1 z^j + c_2 z^k$ for some nonzero scalars $c_1$ and $c_2$.
Then
\[  Z^T(c_1 E^{ij}+ c_2 E^{ik})Z = 2 z^i (z^i)^T \succeq 0,
\]
and
\[  P^T(c_1 E^{ij}+ c_2 E^{ik})Z = (p^i+c_1 p^j + c_2 p^k) (z^i)^T.
\]
Since
$\N( z^i (z^i)^T)$ = $\N((z^i)^T) \subseteq $ $\N( (p^i+ c_1 p^j+c_2 p^k) (z^i)^T)$, it follows
from Equation (\ref{eqj1}) that there exists $t < 0$ such that
$D + t (c_1 E^{ij}+ c_2 E^{ik})$ is an EDM, and thus entries $d_{ij}$ and $d_{ik}$ are jointly yielding.

On the other hand, assume that there exist no nonzero scalars $c_1$ and $c_2$ such that
$z^i= c_1 z^j+ c_2 z^k$, and assume, to the contrary, that entries $d_{ij}$ and $d_{ik}$ are jointly yielding.
Then there exist $t_1 \neq 0$ and $t_2 \neq 0$ such that $D + t_1 E^{ij} + t_2 E^{ik}$ is an EDM.
Thus it follows from Equation (\ref{eqj1}) that
\begin{eqnarray*}
& & - Z^T(t_1 E^{ij}+ t_2 E^{ik})Z  \succeq \bz, \mbox{ and } \\
& &  \N(Z^T(t_1 E^{ij}+ t_2 E^{ik})Z) \subseteq \N(P^T(t_1 E^{ij}+ t_2 E^{ik})Z).
\end{eqnarray*}
But
\[
Z^T(t_1 E^{ij}+ t_2 E^{ik})Z = z^i ( t_1 z^j + t_2 z^k)^T + ( t_1 z^j + t_2 z^k) (z^i)^T
\]
and
\[
P^T(t_1 E^{ij}+ t_2 E^{ik})Z = p^i ( t_1 z^j + t_2 z^k)^T + ( t_1 p^j + t_2 p^k) (z^i)^T.
\]
We have three cases to consider here.

(i) both $z^i$ and $t_1 z^j+ t_2 z^k$ are nonzero. Thus it follows from
Proposition~\ref{prop1} that $Z^T(t_1 E^{ij}+ t_2 E^{ik})Z$ is indefinite, a contradiction.

(ii)  $z^i=\bz$ and $t_1 z^j+ t_2 z^k \neq \bz$. Thus  $Z^T(t_1 E^{ij}+ t_2 E^{ik})Z=\bz$. Moreover, it follows from
Proposition \ref{propz=0} (part 1) that $P^T(t_1 E^{ij}+ t_2 E^{ik})Z$ = $ p^i ( t_1 z^j + t_2 z^k)^T \neq \bz$. Thus
$\N(Z^T(t_1 E^{ij}+ t_2 E^{ik})Z) \not \subseteq \N(P^T(t_1 E^{ij}+ t_2 E^{ik})Z)$, also a contraction.

(iii) $z^i \neq \bz$ and $t_1 z^j+ t_2 z^k = \bz$. Then $z^j=-t_2 z^k/ t_1$. Hence, it follows from
Proposition \ref{propz=0} that $p^j+ t_2 p^k/ t_1 \neq \bz$. Thus, $P^T(t_1 E^{ij}+ t_2 E^{ik})Z$ =
$( t_1 p^j + t_2 p^k) (z^i)^T \neq \bz$. Then the same argument in case (ii) leads to
a contradiction.

Therefore, $d_{ij}$ and $d_{kl}$ are not jointly yielding.
\epr

The following two corollaries follow immediately if $D$ is in general position. Recall that the notion of
jointly yielding (or jointly unyielding) is defined only for two unyielding entries in the same row (column).

\begin{cor} \label{corjgp1}
Let $D$ be an $n \times n$ EDM of embedding dimension $r=n-3$. If $D$ is in general position in $\Rs^r$, then
every entry of $D$ is unyielding and every two entries of $D$, in the same column (row), are jointly yielding.
\end{cor}

\bpr
The fact that every entry of $D$ is unyielding follows from Corollary \ref{corgp2}.
In this case, $z^1, \ldots, z^n$ are in $\Rs^2$ since $(n-r-1)=2$. Let
$i, j$, and $k$ be any three distinct indices in $\{1,\ldots,n\}$.
Since $D$ is in general position,
it follows from Lemma \ref{lemgp} that any 2 of $z^1, \ldots, z^n$ are linearly independent.
In particular,
$\{z^j, z^k\}$ is a basis of $\Rs^2$ and thus
$z^i= c_1 z^j + c_2 z^k$ for some nonzero scalars $c_1$ and $c_2$. Therefore,
entries $d_{ij}$ and $d_{ik}$ are jointly yielding.
\epr

\begin{cor} \label{corjgp2}
Let $D$ be an $n \times n$ EDM of embedding dimension $r \leq n-4$. If $D$ is in general position in $\Rs^r$,
then every entry of $D$ is unyielding and every two entries of $D$, in the same column (row), are  jointly unyielding.
\end{cor}

\bpr
The fact that every entry of $D$ is unyielding follows from Corollary \ref{corgp2}. Now
assume, to the contrary, that two entries of $D$, say $d_{ij}$ and $d_{ik}$, are jointly yielding.
Thus, it follows from Theorem \ref{thmchar2} that $z^i=c_1 z^j + c_2 z^k$ for some nonzero scalars $c_1$ and $c_2$.
Note that in this case, $(n-r-1) \geq 3$.
Hence, any $(n-r-1) \times (n-r-1)$ submatrix of $Z$ containing $(z^i)^T$, $(z^j)^T$ and $(z^k)^T$ is singular. This
contradicts  Lemma \ref{lemgp} and the proof is complete.
\epr

In the following two theorems, we provide upper and lower bounds for each jointly yielding pair of entries
of an EDM $D$.

\begin{thm}\label{thmjiz=0}
Let $D$ be an $n \times n$, $n \geq 4$, EDM of embedding dimension $r \leq n-2$ and let $Z$ and $P$ be a
Gale matrix and a configuration matrix of $D$ respectively, where $P^Te=\bz$.
Further, Let $S= P (P^TP)^{-1}$ and let $s^i$ be the $i$th row of $S$, i.e., $s^i = (P^TP)^{-1} p^i$.
If $z^i = \bz$, $z^j \neq \bz$, $z^k \neq \bz$ and $c_1 z^j + c_2 z^k = \bz$ for some
nonzero scalars $c_1$ and $c_2$. Then   $D + t (c_1 E^{ij} + c_2 E^{ik}$) is an EDM if and only if
\[
\frac{2}{\lambda_r} \leq t \leq \frac{2}{\lambda_1},
\]
where $\lambda_1 = (s^i)^T (c_1 s^j + c_2 s^k) + || s^i|| \; ||c_1 s^j + c_2 s^k||$ and
$\lambda_r = (s^i)^T (c_1 s^j + c_2 s^k) - || s^i|| \; ||c_1 s^j + c_2 s^k||$.
\end{thm}

\bpr The proof is analogous to that of Theorem \ref{thmizz=0}. In particular,
$D + t (c_1 E^{ij} + c_2 E^{ik})$ is an EDM iff
\[
2(P^TP)^2 - t ( p^i (c_1 p^j + c_2 p^k)^T + (c_1 p^j + c_2 p^k) (p^i)^T) \succeq \bz,
\]
iff
\beq \label{eqjz=0}
2 I_r - t ( s^i (c_1 s^j + c_2 s^k)^T + (c_1 s^j + c_2 s^k) (s^i)^T) \succeq \bz.
\eeq
Now it follows from Proposition \ref{propz=0} (parts 1 and 3) that $p^i \neq \bz$ and $c_1 p^i + c_2 p^j \neq \bz$
since $z^j=-c_2 z^k/ c_1$. Also, it follows from  Proposition \ref{propz=0} (part 2) that
$p^i$ is not parallel to $c_1 p^j + c_2 p^k$. Therefore, it follows from Proposition \ref{prop1} that
Equation (\ref{eqjz=0}) holds iff $2 - t \lambda_1 \geq 0$ and $2 - t \lambda_r$ where
$\lambda_1 > 0 > \lambda_r$ are the nonzero eigenvalues of $s^i (c_1 s^j + c_2 s^k)^T + (c_1 s^j + c_2 s^k) (s^i)^T$.
\epr

Note that the result of Theorem \ref{thmjiz=0} reduces to that of Theorem \ref{thmizz=0} if we set $c_2=0$.
\begin{exa}
Consider the EDM  $D = \left[ \begin{array}{cccc} 0 & 17 & 16 & 17 \\ 17 & 0 & 1 & 4 \\
                                                   16 & 1 & 0 & 1   \\ 17 & 4 & 1 & 0
                                                   \end{array} \right] $ of embedding dimension~$2$.
A configuration matrix and a Gale matrix  of $D$ are
\[ P=  \left[ \begin{array}{rr} 0 & 3 \\ -1 & -1 \\ 0 & -1 \\ 1 & -1 \end{array} \right] \mbox{ and }
 Z= \left[ \begin{array}{r} 0 \\ 1  \\ -2 \\  1 \end{array} \right].
\]
Moreover,
\[
S= P (P^TP)^{-1} = \frac{1}{12}\left[ \begin{array}{rr} 0 & 3 \\ -6 & -1 \\ 0 & -1 \\ 6 & -1 \end{array} \right].
\]
Then $d_{12}$ and $d_{13}$ are unyielding. However, since $z^1=0= 2 z^2 + z^3$, it follows that
 $d_{12}$ and $d_{13}$ are jointly yielding. In this case,
$\lambda_1= (s^1)^T (2 s^2+s^3) + ||s^1|| \; ||2s^2+s^3|| = (\sqrt{17} - 1)/16$ and $\lambda_2= - (\sqrt{17}+1)/16$. Hence,
$D+ t (2 E^{12} + E^{13})$ is an EDM iff $-2 (\sqrt{17} -1) \leq t \leq 2 (\sqrt{17} + 1)$.
\end{exa}

\begin{thm}\label{thmjizpz}
Let $D$ be an $n \times n$, $n \geq 4$, EDM of embedding dimension $r \leq n-2$ and let $Z$ and $P$ be a
Gale matrix and a configuration matrix of $D$ respectively, where $P^Te=\bz$.
Further, Let $S= P (P^TP)^{-1}$ and let $s^i$ be the $i$th row of $S$, i.e., $s^i = (P^TP)^{-1} p^i$.
If $z^i \neq 0$, $z^j \neq 0$, $z^k \neq 0$ and $z^i= c_1 z^j + c_2 z^k$ for some
nonzero scalars $c_1$ and $c_2$. Then   $D + t (c_1 E^{ij} + c_2 E^{ik})$ is an EDM if and only if
\[
\frac{-4}{||s^i - c_1 s^j - c_2 s^k||^2} \leq t \leq 0.
\]
\end{thm}

\bpr
The proof is analogous to that of Theorem \ref{thmizpz}. In particular,
$D + t (c_1 E^{ij} + c_2 E^{ik})$ is an EDM iff
\[
 t \leq 0 \mbox{ and } 2 (P^TP)^2 + \frac{t}{2} (p^i-c_1 p^j - c_2 p^k) (p^i-c_1 p^j - c_2 p^k)^T \succeq \bz,
\]
which is equivalent to
\[
 t \leq 0 \mbox{ and } 2 I_r + \frac{t}{2} (s^i-c_1 s^j - c_2 s^k) (s^i-c_1 s^j - c_2 s^k)^T \succeq \bz.
\]
Moreover, it follows from Proposition \ref{propz=0} (part 4)  that $s^i-c_1 s^j - c_2 s^k \neq \bz$,
and thus the proof is complete.
\epr

Also note that the result of Theorem \ref{thmjizpz} reduces to that of Theorem \ref{thmizpz} if we set $c_2=0$.

\begin{exa}
Let $D$ be the EDM considered in Example \ref{exalast}.
Then
\[
S= P (P^TP)^{-1} = \frac{1}{26}\left[ \begin{array}{r} -7 \\ -2  \\ 3 \\ 8 \\ -2 \end{array} \right].
\]
Then $d_{13}$ and $d_{14}$ are unyielding. However, since $z^1=  z^3 + 2 z^4$, it follows that $d_{13}$ and
$d_{14}$ are jointly yielding. In this case,
$s^1 - s^3 - 2 s^4 = -1$. Hence,
$D+ t ( E^{13} + 2 E^{14})$ is an EDM iff $-4 \leq t \leq 0 $.
\end{exa}



\begin{thebibliography}{1}

\bibitem{alf16}
A.~Y. Alfakih.
\newblock Universal rigidity of bar frameworks via the geometry of
  spectrahedra.
\newblock  {\em J. Global Optim.}, 67:909--924, 2017.

\bibitem{alf07}
A.~Y. Alfakih.
\newblock On dimensional rigidity of bar-and-joint frameworks.
\newblock {\em Discrete Appl. Math.}, 155:1244--1253, 2007.

\bibitem{cri88}
F.~Critchley.
\newblock On certain linear mappings between inner-product and squared distance
  matrices.
\newblock {\em Linear Algebra Appl.}, 105:91--107, 1988.

\bibitem{gal56}
D.~Gale.
\newblock Neighboring vertices on a convex polyhedron.
\newblock In {\em Linear inequalities and related system}, pages 255--263.
  Princeton University Press, 1956.

\bibitem{gow85}
J.~C. Gower.
\newblock Properties of {E}uclidean and non-{E}uclidean distance matrices.
\newblock {\em Linear Algebra Appl.}, 67:81--97, 1985.

\bibitem{gru67}
B.~Gr{\"{u}}nbaum.
\newblock {\em Convex polytopes}.
\newblock John Wiley \& Sons, 1967.

\bibitem{sch35}
I.~J. Schoenberg.
\newblock Remarks to {M}aurice {F}r\'{e}chet's article: Sur la d\'{e}finition
  axiomatique d'une classe d'espaces vectoriels distanci\'{e}s applicables
  vectoriellement sur l'espace de {H}ilbert.
\newblock {\em Ann. Math.}, 36:724--732, 1935.

\bibitem{wsv00}
H.~Wolkowicz and R.~Saigal and L.~Vandenberghe.
\newblock Handbook of Semidefinite Programming. Theory, Algorithms and Applications.
\newblock Kluwer, 2000.

\bibitem{yh38}
G.~Young and A.~S. Householder.
\newblock Discussion of a set of points in terms of their mutual distances.
\newblock {\em Psychometrika}, 3:19--22, 1938.

\end{thebibliography}

\end{document}